# Cyber-Physical Microgrids: Toward Future Resilient Communities

Tuyen V. Vu, *Member, IEEE,* Bang H. L. Nguyen, *Student Member, IEEE,* Zheyuan Cheng, *Student Member, IEEE,* Mo-Yuen Chow, *Fellow, IEEE,* Bin Zhang, *Senior Member, IEEE*

*Abstract*—Microgrids can be isolated from large-scale power transmission/distribution systems (macrogrids) to deliver energy to their local communities using local energy resources and distribution systems when power outages occur in the macrogrids. In such situations, microgrids could be considered as the last available resource to provide energy to critical infrastructure. Research in monitoring and control of microgrids has been ongoing for the last two decades to protect and enhance communities' socio-economic performance. However, of increasing concern are the possible cyber-physical threats that could disrupt the provision of macrogrids' energy services to critical infrastructure and consequently impact the resilience and sustainability of communities. As cyber-physical systems, microgrids are not immune to these threats. Advanced monitoring and control are critical for real-time operations of microgrids and, therefore, directly influence communities' resilience. Research trends in monitoring have recently shifted from normal situational awareness in forecasting, state estimation, and prediction to anomalies' analysis and cyber-physical attacks' detection to support resilient control systems. In addition, confounding the interpretation of research findings is the lack of a widely accepted definition, analytical methods, and metrics to consistently describe the resilience of power grids, especially for microgrids. This paper provides an overview of current research in microgrid resilience and presents an outlook for future trends.

## I. Overview of Cyber-Physical Microgrids

**MICROGRIDS play a significant role in ensuring resilient and sustainable energy to communities:** Critical community infrastructures, such as hospitals, water treatment plants, and emergency and military services, rely heavily on power and energy services for their resilience and survivability, especially in the face of natural disasters. Catastrophic events such as hurricanes, blizzards, thunderstorms, and earthquakes can severely damage community infrastructures and result in power outages (among other consequences) that may take weeks to resolve. The economic and social hardships [1] of power disruptions caused by natural disasters in the U.S. is significant and has been estimated to cost $25–70 billion annually [2]. An emerging strategy to mitigate the consequences of power outages and promote resilience of the power grid is through the introduction of microgrids. As controllable entities, microgrids can be seamlessly connected or disconnected from macrogrids once

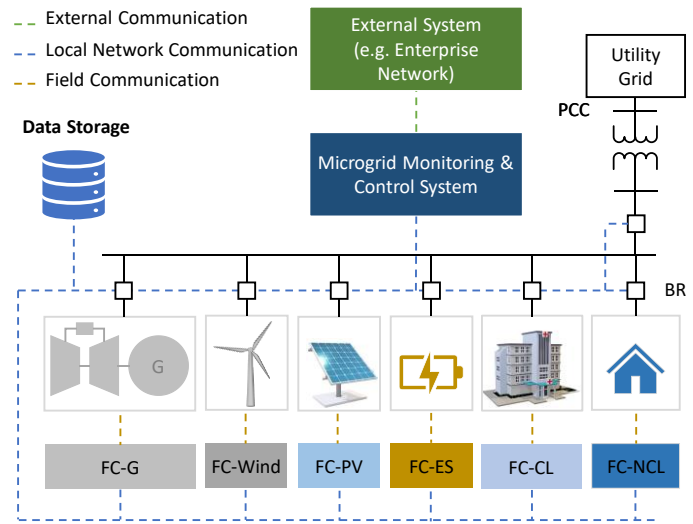

FIGURE 1 - Cyber-physical Microgrid. Abbreviations: BR: breaker, CL: critical load, ES: energy storage, FC: field controller, NCL: noncritical load, PCC: point of common coupling, PV: photovoltaic.

outage events occur and so maintain required energy services to meet local demands. GTM Research forecasts that microgrid capacity in the U.S. will grow from 3.2 GW in 2017 to 6.5 GW in 2022, a 14.1 percent compound annual growth rate [3]. Worldwide, as of the second quarter of 2019, Navigant Research identified nearly 4500 microgrid projects (representing almost 27 GW power capacity) that have been or will be installed [4].

**Advanced monitoring and control are required for optimal operations of microgrids:** Microgrids are cyber-physical systems (CPS) that contain a wide variety of interconnected devices that measure, control, and actuate distributed energy resources (DER), loads, and power distribution devices. A typical microgrid is shown in FIGURE 1, where the monitoring and control system communicates internally with local devices via a local area network and externally with its enterprise network or power distribution systems via a wide area network. Since microgrids can be considered as an integrated power and energy system with DER, loads, and distribution automation devices, their monitoring and control can be as complex as in a bulk transmission system. Monitoring functions provide overarching information about the current system states; this information is also used to predict the system's future states and events. Based upon the nature of the information monitored, a properly performing control system would provide optimal control actions, ensuring the resilience of the microgrid and critical community infrastructures.

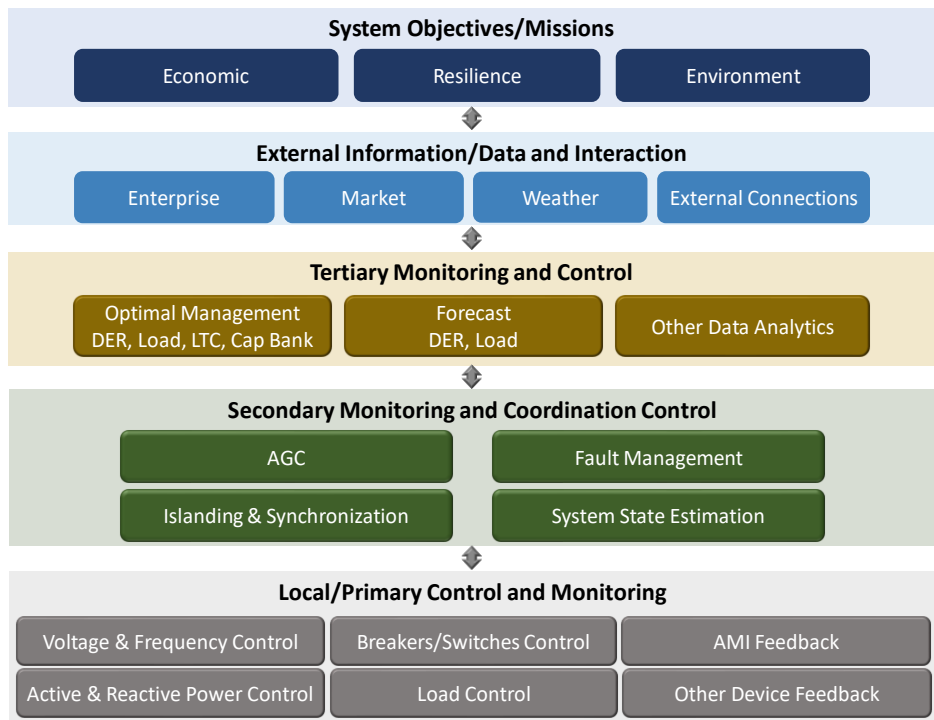

FIGURE 2 –Monitoring and control functionalities in microgrids.

Recent reviews have covered forecasting and restoration methods for generic power systems' resilience under natural disasters [5], networked microgrids in enhancing power system resilience against extreme events [6], or highlighted specific microgrid characteristics and surveyed research directions and challenges [7]. However, to the best of our knowledge, a comprehensive review of work on critical monitoring and control functions that consider future trends pertaining to the resilience of cyber-physical microgrids has not yet been published.

The rest of the paper is organized as follows: In Section 2, we (1) present an overarching multilayer architecture that covers microgrids' control and monitoring functionalities; (2) review critical monitoring functions including forecasting, system state estimation, and anomaly detection; and (3) summarize recent advanced control efforts that contribute to the resilience of a three-layer control system. Informed by this review, in Section 3 we identify future research directions to promote resilient microgrids against emerging cyber-physical threats. In Section 4, we summarize the review and research directions.

## II. Recent Trends in Monitoring and Control for Cyber-Physical Microgrids

### A. Monitoring and Control Architectures

An architecture that lays out monitoring and control functions for microgrids is critical as it provides researchers and engineers a framework for effective research and development efforts. The recent IEEE standard 2030.7-2017 identifies critical control functions for a centralized microgrid controller [8]. However, functions for more general purposes (i.e., those not specific to either centralized or distributed control) that are either consistent with or complementary to this standard are illustrated in FIGURE 2. This architecture has three main layers of monitoring and control. The first is the primary control and monitoring layer, where basic functions occur: DER (active/reactive power, frequency/voltage control), breakers/switches (on/off and protection functions), loads' control (curtailment functions), advanced metering infrastructure (AMI) devices (voltage, current, and power measurement), and feedback from other devices (DER voltage, current, and power)). In the secondary layer, information feedback from the primary level can be processed for system state estimation to monitor the system's behavior for control purposes. Control functions could also be performed including the frequency-restoration function (automatic generation control (AGC)), islanding operations and synchronization, and fault management (coordination of breakers/switches). In the tertiary layer, advanced monitoring functions such as forecasting (DER and load forecasting) can be performed, and other data analytics such as anomaly detection of AMI information (e.g., smart meters) can be conducted. Based upon information obtained from monitoring functions, the optimal management strategies for DER, load, load tap changer (LTC), and capacitor bank could be performed to achieve the system's objectives or missions, whether they pertain to economic, resilience, or environmental attributes. The socio-economic performance measures associated with these objectives would be specified by the microgrid operator.

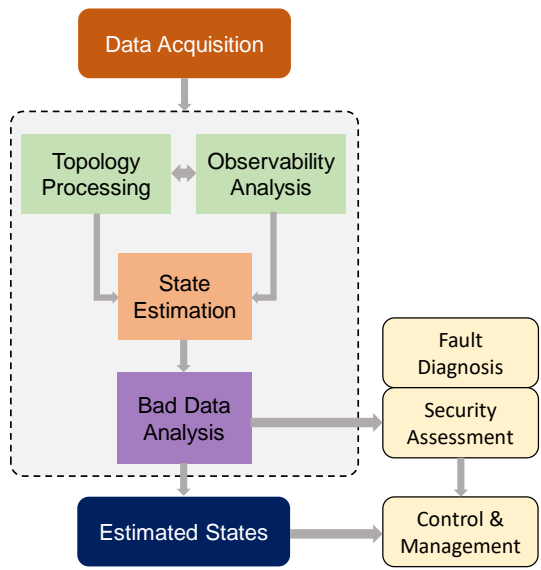

FIGURE 3 – System state estimation.

TABLE 1. SE FRAMEWORKS AND TECHNIQUES MERITS AND DRAWBACKS.
NOTE: (+) INDICATES ADVANTAGES AND (-) INDICATES DISADVANTAGES.

| SE Framework | Robustness | Nonlinearity |
|---|---|---|
| **SSE**: (+) Simple, robust to uncertainties, (-) Do not track system dynamics. **TSE**: (+) Simple, suitable for stiff system states, (-) Cannot track variations of flexible loads and generation units, and low performance with weak grids. **FASE**: (+) involving states/load forecasting to SE, good performance at smooth change, (-) Neglect dynamics, inaccurate results under highly state fluctuations. **DSE**: (+) Provide best performance with state space representations, (-) Complex and high computation, where dynamic models are not always available. | **WLS**: (+) simple, popular for static SE, (-) sensitive to outliers. **KF**: (+) same cost function of WLS, popular for dynamic SE, (-) sensitive to outliers. **LAV**: (+) Robust against outliers, less sensitive to parameter errors (-) high computing cost, sensitive to attacks and measurement errors. **GML**:(+) robust against outliers, (-) sensitive to parameter errors. **H-infinity**: (+) limit system uncertainties (-) lack robustness to outliers and non-Gaussian noise | **Linearization**: Taylor approximation: (+) Constant Jacobian matrix with rational state variables, relatively simple (-) Low performance with highly nonlinear systems. **Noises Propagation: Unscented, ensemble, particle transformations**: (+) High performance with highly nonlinear systems, (-) Relatively complex. |

## B. Monitoring Systems

The monitoring functionalities in the three possible layers provide situational awareness to support control decisions. In the following, we review trends in monitoring functions, specifically state estimation at the secondary level, DER and load forecasting at the tertiary level, and anomaly detection at any monitoring level are critical for overall system situational awareness.

### 1) System State Estimation

State estimation (SE) is important for microgrid operations, particularly for energy management and secondary control of system voltage and frequency. State estimation usually processes noisy measurements to extract accurate states of power networks [9]. FIGURE 3 shows SE in relation to other functions in power systems. As a small power network, a microgrid requires SE, specifically distribution system state estimation (DSSE) algorithms, to extract its true states. Unlike conventional SE applied to transmission systems, DSSE algorithms suffer from the following limitations [10]: (1) Its observability is low due to the massive number of nodes, resulting in thousands of states to be estimated versus an insufficient number of measurements [10]. (2) Highly dynamic load profiles result in inaccurate forecasted load (pseudo measurements) [11]. And (3), DSSE is unable to decouple active and reactive powers in their formulations because of three-phase unbalanced operations and low reactance/resistance ratios [12]. In the literature, two types of DSSE have been reported, categorized as model-based or data-driven [13].

**Model-based methods** rely on mathematical models of power networks and measurements. Conventionally, a weighted least-squares (WLS) problem can be formulated for static SE [14]. However, the WLS process needs to reinitialize at every time-step with new measurements making future estimated states independent from historical states. Therefore, this method detects bad data only within the snapshot measured. Consequently, it is vulnerable to false data injection attacks (FDIA) and other malicious data manipulation activities [15]. As an alternative, dynamic SE (DSE) algorithms have been studied to provide better situational awareness as they involve the state-space model of power networks, revealing the system states' evolution [16]. Various DSE methods have been applied to estimate system states [17]; among these are Kalman filter (KF) variations, such as extended KF [18], ensemble KF [19], unscented KF [20], and particle filter [21], are frequently applied. Dynamic SE can track state changes and detect bad data using a normalized innovation or chi-square test [20]. However, DSE works based on accurate dynamic power system models, which are not always available. Therefore, DSE can be applied to synchronous generators' state estimates, whereas the tracking SE (TSE) [16] is applied to network states provided that their operation is quasi-steady [22]. For loads that are forecasted, the forecasting-aided SE (FASE) [23] is applied. To improve the robustness of DSE under the presence of outliers, least-absolute-value (LAV) and generalized maximum-likelihood (GML) methods have been proposed [22], [18]. H-infinity based filtering has also been introduced to bound the system's uncertainties [24]. In DSE, both node voltage and branch current states have been adopted in polar or rectangular form. The rectangular form yields a simpler formulation and more effective computation than the polar form does [25]. State estimation also can be realized using centralized or decentralized approaches. In decentralized approaches, the whole network is partitioned into subareas, where local estimators can be executed in parallel [26]. Centralized approaches have higher computation and communication burdens than decentralized methods [27] in part because measurement redundancy is lessened by partitioning of the whole power network into subnetworks, and data fusion algorithms are required to reconstruct the global estimated states [26]. The state estimation frameworks with nonlinearity and robustness properties are compared as shown in Table 1. Based on the comparison, future work can focus on the enhancement of the accuracy and robustness of DSE, especially under Gaussian and/or nonGaussian noise, cyber-attacks, and bad data [24]. Research in distributed schemes for microgrids also needs to keep pace with recent developments [28]. In addition, future DSE should cover fault diagnosis and security assessment and support fault tolerant control to enhance system resilience [29], [30].

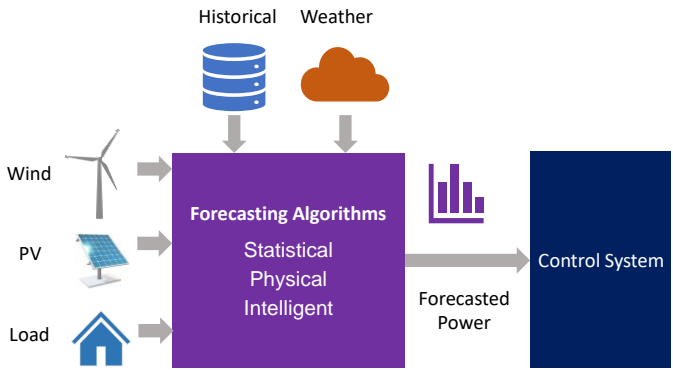

FIGURE 4 – Forecasting in microgrids.

**Data-driven approaches** have recently become an attractive research direction since they require no physical models to estimate system states. The major trend in these approaches is to utilize artificial neural networks (ANN) to represent either parts of or the entire power network through training processes using extensive data collected from various resources. These ANN-based methods extract power networks' operational patterns and represent them as weights of neural nodes. This approach can improve the robustness of the SE against erroneous measurements with either a probabilistic neural network [31] or a parallel distributed processing model [32]. This approach when combined with traditional model-based methods results in improved observability with enhanced pseudo-measurement generation [33], [34]. Deep neural networks such as the stacked auto-encoder have also been applied recently to provide AC state estimation against cyber-attacks [35]. Besides ANN, K-nearest neighbors search, supervised learning, and kernel trick methods have also been used for grid current state inference [36]. The main advantage of data-driven approaches over traditional SE (model-based) methods is that data-drive approaches do not depend on the system model, which could change during system operations. Therefore, they scale more readily than do model-based approaches. However, model-based approaches tend to be more accurate than data-driven based approaches. Despite the promise motivating the emerging usage of data-driven approaches, they require further research to prove their feasibility, applicability, and superiority compared to traditional methods.

*2) Distributed Energy Resources and Load Forecasting*

The DER (particularly PV and wind) and load forecasting (FIGURE 4) are critical for optimal operations of microgrids. In the following, we review recent literature on forecasting algorithms for loads and DER.

**Load forecasting.** Load forecasting (LF) is classified as either long-term (LTLF) or short-term (STLF). Long-term LF focuses on load operations over a range of weeks or more and is used for resource planning; STLF focuses on load operations over a range of hours [37] and is used for real-time microgrid optimization, especially for the energy management system (e.g., peak shaving and demand response). However, STLF in microgrids is challenging because of the high variability and nonlinearity of load demands compared to bulk power systems [38]. Load forecasting can be performed using either statistical or intelligent approaches [39]. One popular statistical method, regression-based techniques, has been applied to probabilistic load forecasting using quantile regression averaging [40]. Other statistical methods include time-series techniques such as autoregressive moving average (ARMA) [41] or autoregressive integrated moving average (ARIMA) [42]. However, time-series statistical methods are infrequently used because variations in microgrid load demand are high compared to bulk power systems. Besides time-series and other statistical approaches, intelligent approaches such as support vector machine (SVM) have been reported for load forecasting of buildings [43]. Although the use of SVM for load forecasting may be applicable to microgrids, this has not yet been reported. However, the utilization of ANN for demand forecasting in microgrids has increased. Examples of ANN applications include their use in a bilevel forecasting structure, a combined neural network, an evolutionary algorithm, and a differential evolution algorithm [38], deep neural networks with multilayer perceptron [44], combined multilayer perceptron models [45], and a self-recurrent wavelet neural network [39].

**DER forecasting.** Forecasting the intermittent energy output of DER, including solar PV and wind, is as critical as load forecasting for optimal microgrid operations. However, because of relatively small geographical space and power size available for prediction (and so sensitivity to changes in the local environment) forecasting energy output from PV and wind is more challenging than in bulk power systems. Forecasting durations of wind and solar energy range from very short-term (intra-hour), short-term (intra-day), or long-term (day-ahead) [46], [47]. Common approaches for solar or wind forecasting are statistical, intelligent, physical-based, or hybrid approaches [48], [49].

In solar forecasting, the appropriateness of each forecasting approach depends upon the spatial-temporal relationship of the application; statistical and intelligent approaches are more applicable to microgrids while physical-based approaches based on numeric weather prediction (NWP) are more applicable to bulk power systems. Statistical approaches include persistence and regression-based methods (ARMA and ARIMA) in which cloud images and satellite data are used [50].

Intelligent approaches include SVM and ANN [51], [52]. Recent literature also shows a trend toward utilizing deep ANN for solar forecasting to improve prediction accuracy in microgrids [53]. In contrast, methods that use NWP models, such as European Centre for Medium-Range Weather Forecasts (ECMWF), Fifth-Generation Penn State/NCAR Mesoscale Model (MM5), and WRF, are more applicable to wide-area solar forecasting having resolution of tens of kilometers [46]; therefore, they might not be applicable to microgrids.

Statistical methods such as persistence, ARMA, and ARIMA are also used to forecast wind-farm energy output, although their time horizon is normally limited to intra-hour and intra-day forecasts [54]. Intelligent methods including SVM and ANN have also been investigated for improved forecasting accuracy. The prediction accuracy of SVM has been shown to be significantly better than the persistence model's [55]. Over the last decade, ANN has been extensively investigated to improve the accuracy of wind-energy prediction using architectures with many layers and variations [56], [57]. Wind forecasting using physical-based approaches relies on NWP data to provide atmospheric variation over time, which becomes the input for the physical wind model to predict energy output [58]. Similar to solar forecasting, NWP data can be generated from ECMWF, MM5, or WRF models; however, the low space-resolution makes this physical-based method of limited applicability to microgrids. Most recent papers focus on forecasting for large wind systems (windfarms); microgrids, in contrast, have more distributed wind with relatively small energy output. Therefore, the applicability of these methods to microgrids needs further investigations.

*3) Anomaly and Intrusion Detection*

Microgrids, as complex cyber-physical systems (CPS), comprise physical power networks and information and communication technology (ICT) infrastructures that support control and monitoring [59]. In such CPS, detection of anomalies and intrusion attacks is the principal strategy to ensure secure microgrid operations [60]. Anomaly detection should be considered at all levels of monitoring. Anomaly detection approaches for cyber-physical security of power systems can be classified into two categories: physics-based and cyber-based [61].

**Physics-based methods** build indicators of faults, attacks, and anomalies using data related to the physical power network [62]. Physical data received from sensors are processed using the SE tools. The differences between estimated and received data, called residuals, are derived through bad data identification algorithms [60]. These algorithms are statistical tests pertaining to two main categories: stateless and stateful. In a stateless test, a metric indicator, built on a residual or innovation vector at a single time step, is compared to a threshold for detection. Traditionally, in static SE, measurement errors are evaluated by a residual vector ($J(\hat{x})$ detector) [14] or the largest normalized residual indicator (LNR) as $l_\infty$-norm [63]. Alternatively, the generalized likelihood ratio detector has been proved to have better performance than the traditional ones on large sample sizes [64]. In dynamic SE, measurement errors and sudden changes are normally detected via the normalized innovation ratio (NIR) [65]. However, NIR does not involve bad data in states or inputs; therefore, NIR may be bypassed by these attacks. By considering accumulated effects over time, a Chi-square detector can be used to detect soft failures such as instrument bias shift [66]. However, as a Chi-square detector can be bypassed by false data injection attacks (FDIA), the Euclidean distance metric, which evaluates the deviation of measured and estimated data, is introduced to enhance detection capability [66]. However, since straight-forward stateless tests use only the spatial relationship of data for detection, both the false positive (false alarms) rate is high as is the false negative as stealthy attacks may be missed. Alternatively, a stateful test better indicates anomalies [67] by tracking the historical changes of design metrics over time; in other words, it leverages the temporal relationship of measurements. For instance, to detect FDIA nonparametric cumulative sum (CUSUM) statistics accumulate the expected value of observations [68]. However, under unknown statistic attacks, CUSUM can be bypassed. Applying the generalized likelihood ratio test (GLRT) approach with CUSUM can resolve this problem [69]. Another metric, the signal temporal logic, is applied to detect anomalies [70]; but, this test only compares the measured value with its threshold, which could be fooled by stealthy attacks. To detect FDIA, Absolute and Kullback-Leibler distances are employed to track the dynamics of adjacent measurements [71].

Besides these tests and indicators, data-driven and machine-learning tools have recently been utilized for anomaly detection. Such tools include semisupervised learning [72], deep autoencoders [73], naive Bayes classification [74], gate recurrent unit with multi-layer perceptron [67], principal component analysis [75], advanced multigrained cascade forest algorithm [76], and reinforcement learning [77]. Although this area is attracting many researchers, the efficiency and computational cost of these new methods need to be evaluated and compared with traditional methods to identify the most efficient and cost-effective physics-based strategies for anomaly detection.

**Cyber-based methods** identify anomalies by leveraging IT data extracted from electronic devices and communication channels and can be classified into network-based and host-based approaches [78]. Network-based methods capture and assess communication packets and network behavior, whereas host-based methods identify intrusion footprints within a host-device by evaluating activities' logging, the integrity of system files, and finite machine states [79]. Regarding the host-based approach, an embedded intrusion detection method was proposed in [80] for intelligent electronic devices (IED) in substations that monitor all incoming Generic Object-Oriented Substation Events (GOOSE) and Sampled Values (SV) messages to enhance cybersecurity. Another host-based algorithm has been proposed based on intrusion footprints in user-interfaced computers and IED [81]. Regarding network-based approaches, [82] proposed an intrusion detection system (IDS) for a ZigBee home area network, where network features are analyzed. In [83], a method was proposed that transforms behavior rules for network devices to a state machine and compares their behavior specifications. A hierarchical IDS framework was proposed in [84], which employed an SVM and an artificial immune system to analyze network traffic at every layer of a smart grid. A similar method was also presented in [85] but for individual AMI devices. Besides host- and network-based approaches, hybrid anomaly detection systems have been proposed to cover both the host and network sides for substations [78]. In [86], a SCADA-specific IDS was introduced (based on IEC61850 standard) to cover detection mechanisms based on access control, protocol whitelisting, the cyber model, and combined parameters.

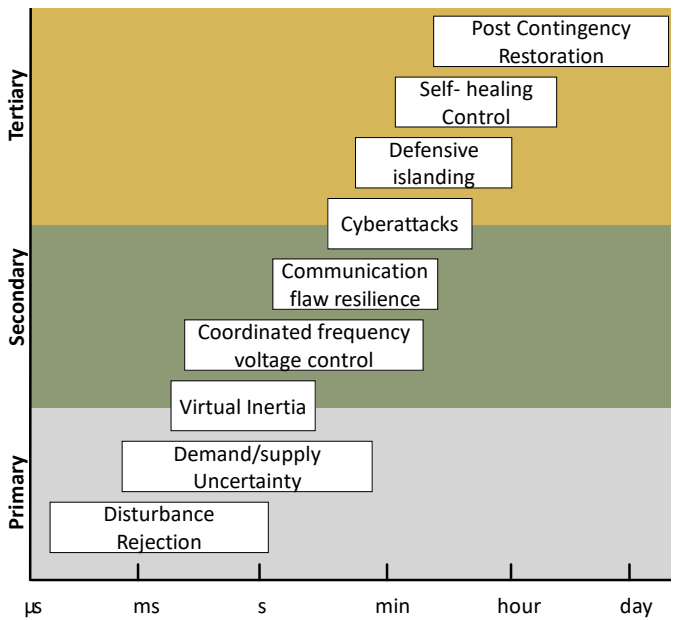

FIGURE 5 – Key resilience control objectives.

To attain a robust and secure microgrid protected from various kinds of cyber-attacks, a defense-in-depth architecture is recommended when designing a microgrid's control system; this strategy focuses on multiple security layers to detect and isolate attacks [87], [59]. Software-defined communication, an emerging networking paradigm, has also been implemented to enhance resilience of smart grids with remarkable results and may also be suitable for microgrids [88], [89].

### C. Control Systems

As noted, a microgrid control system can be divided into three layers: primary, secondary, and tertiary [90]–[92]. This article adopts widely acknowledged definitions for these layers from [93], which is also consistent with the standard P2030.7-2017. The primary layer includes device-level real-time feedback controls based on local measurements. The secondary layer includes system-level voltage and frequency regulations. The tertiary layer consists of slow time-scale system-level controls (e.g., energy management system (EMS), load restoration, and system reconfigurations).

The concept of microgrid resilience is not yet standardized in the power and energy community [94], [95]. Most published work on resilient control focuses on one or a few control objectives within specific layers. The key resilient control objectives are summarized in FIGURE 5. In the following, we systematically discuss recent literature addressing resilience control of three layers in microgrids against commonly reported threats: natural disasters, disturbances and faults, communication flaws, and cyberattacks [95], [6].

*1) Resilience in Primary Control Layer*

As previously mentioned, the primary layer contains DER control and device feedback. With the growth of DER, including renewable energy sources in microgrids, advanced control is increasingly important. The aim of advanced control is to produce accurate and low distorted voltage and current under large system disturbances, such as faults and highly variable power generation and load demand. Many control strategies have been investigated including traditional PID control [96], state feedback control [97], model predictive control [98], and sliding mode control [99]. Most of these approaches, however, are either sensitive to uncertainties and disturbances [100] or lead to difficulties in low-pass filtering or overstress on switching devices [101].

To asymptotically reject disturbances, resonant control (RSC) based on the internal model principle (IMP) was proposed [102]. Various modified RSC schemes, such as proportional (P)+RSC [103], proportional-integral (PI)+RSC [101], multi-RSC (MRSC), and phase compensation RSC [104], were developed to improve transient response as well as stability and resilience against disturbances. Due to its effectiveness, RSC has become a popular controller for grid converters. A modified RSC which is equivalent to a parallel combination of a feedback controller (e.g., PID, deadbeat, or state feedback) and RSC components at all harmonics has attracted extensive attention [105]. Over the past decades, significant advancements in DER control have been achieved including various phase compensation methods [101], selective harmonic RC [106], odd harmonic RC [107], discrete Fourier transform (DFT)-based RC [106] $6k\pm1$ RC, and more general $nk \pm m$ RC [108].

Although these traditional schemes have demonstrated success, they lack frequency adaptability to accommodate frequency fluctuation from distributed generation resources caused by external disturbances or internal faults. To overcome such limitations, fractional order control [109] and variable sampling/switching period techniques [110] were developed. More recently, universal fractional-order design and software-based virtual variable sampling schemes have been developed, which have greatly improved the reliability and resilience of the primary control layer's operation under frequency fluctuation and various disturbances [111], [112].

Table 2 shows the advantages and disadvantages of the above-mentioned methods in terms of structure, optimality, periodic signal tracking, and frequency fluctuation.

TABLE 2. COMPARISON AMONG CONTROL METHODS FOR PRIMARY CONTROL LAYER.

| Methods | Pros. | Cons. |
|---|---|---|
| PID | Simple structure; easy tuning. | Control result is accurate; requires experience in parameter tuning. |
| State feedback | Simple structure; easy to make system stable. | Cannot achieve zero tracking error theoretically. |
| RSC | High magnitude gain at the desired frequency. | Complicated structure to compensate multiple high order harmonic frequencies. |
| MPC | High accuracy; optimal solution for a finite horizon. | High computation cost; Not take advantage of the reference signal's periodic property |
| RC | Track any periodic signal with known frequency; low harmonic distortion. | Needs to know the accurate frequency; cannot respond to frequency fluctuation. |
| Variable sampling RC | Track any periodic signal with known frequency and frequency fluctuation. | More computation. |

*2) Resilience in Secondary Control Layer*

To improve resilience in the secondary control layer against natural disasters, researchers have focused on the development of advanced control methodologies capable of maintaining and restoring system voltage and frequency during and after large system disturbances such as faults and contingencies. Solutions are described as either centralized or distributed. Centralized solutions reported in the literature include adaptive virtual inertia to stabilize system frequency against high renewable energy penetration [113] and frequency response rescheduling for distributed generators after disturbances [114]. Although these centralized solutions are effective and easy to implement, they are typically vulnerable to single-point-of-failure and do not scale well. To overcome these challenges, distributed control techniques are normally utilized, e.g. consensus-based methods, to coordinate DER to control the system's voltage and frequency [115], [116].

Researchers have also investigated resilience to communication flaws (e.g., packet loss, delay, and link failure) in the cyber layer of microgrids. Both centralized and distributed secondary control paradigms rely on the communication of local nodes' voltage and frequency to stabilize the system's voltage and frequency with accurate power-sharing among DER. One extensively reported centralized approach is based on software-defined network (SDN) [117], [118]; in it, a separate centralized controller is deployed to optimally control the cyber layer to ensure network QoS (quality of service). Other researchers have leveraged distributed peer-to-peer (P2P) communication protocols and noise-resilient state-observer techniques to achieve resilience against communication channel noise and disturbances [119], [120]. However, the scalability of the P2P configuration is largely constrained by the existing communication infrastructure. A hierarchical structure that is locally distributed and globally centralized is viewed as one of the most promising trends in this area.

In response to another growing threat—cyberattacks—various efforts have been made to secure both control information and decisions that are exposed in the communication network. Some researchers approach this topic from a conventional cybersecurity perspective: proposing resilient control solutions to traditional cyberattacks (e.g., denial-of-service attack [121], [122].) This body of work assumes a centralized attack detection and mitigation platform. Alternatively, other researchers have focused on insider attacks that target the system's secondary control functions (e.g., false data injection and malicious interruption of control algorithm.) Distributed consensus control frameworks based on state observers have been proposed to detect and mitigate these attacks [120], [123]. Generally, this field of research typically adopts a bottom-up approach, where researchers propose solutions to specific vulnerabilities in various components in the microgrid control system. A well-recognized cyberattack resilience framework has yet to be developed to systematically integrate and coordinate the proposed countermeasures.

*3) Resilience in Tertiary Control Layer*

Considerable research on the tertiary control layer has been conducted to address natural disaster challenges. The objective of the tertiary layer is to optimally coordinate microgrid's DER to minimize a disaster's impact and to promote efficient recovery from it. Related research on tertiary control focuses defensive and intentional islanding [124], improving survivability and robustness [125], [126] reconfiguring and self-healing [127], [128], and restoring service [129], [130]. These efforts assume a centralized controller with full observability to optimally dispatch the controllable DER and microgrid switches to minimize loss of load under extreme events. Some distributed tertiary control approaches assume a multi-agent system that can cooperatively solve the underlining optimization problem and execute control commands [131], [132]. However, the distributed methods that have been adopted in the field typically handle standard optimization problems (e.g., convex problems). As the complexity of the optimal control problems grows, centralized evolutionary algorithms must be regularly adopted.

Another pressing challenge in the tertiary control layer is to improve system resilience against communication flaws. This is a critical requirement for the distributed, controlled microgrid; contingencies in the cyber layer could result in wrong or infeasible control decisions and algorithm interruption/failure. Centralized approaches typically adopt the previously mentioned SDN technique to optimally route the network traffic and reconfigure communication links to provide acceptable QoS under cyber contingency [88], [133]. Distributed approaches leverage robust distributed optimization solvers, such as consensus-based subgradient methods, to mitigate the negative impact of communication flaws [134].

Recently, numerous papers have presented centralized methods to handle conventional cyber threats (e.g., denial-of-service attacks [121], sniffing attacks [135]), that target the tertiary control communication interface and traffic. In addition to these conventional cyber threats, some recent research focuses on advanced attacks such as false data injection attacks by using the sliding mode observer-based technique [136], an unscented Kalman filter [137], and a reputation-based system [138]–[140]. Finally, since tertiary control is closely tied to the economic objectives of for-profit entities or the missions of other entities, malicious attacks that specifically target energy markets have been discussed in the literature [141]–[143].

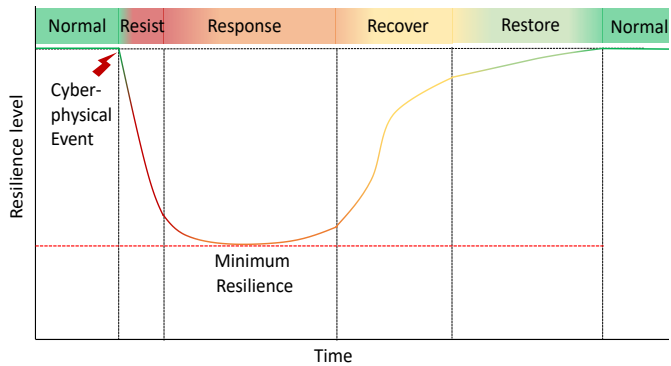

FIGURE 6 – Resilient curve.

### III. FUTURE DIRECTIONS

In the following, we provide our outlook on future research directions that would allow for resilient microgrids against cyber-physical threats.

**Emerging cyber-physical threats:** While most power outages are the result of extreme natural events, there is increasing concern about outages caused by man-made incidents including cyber-physical attacks. The Ukraine power grid attack in December 2015 that left about 230 thousand people without electricity for several hours is a now-famous case of a cyber-attack in power systems [144]. Consequently, an extreme outage caused by a natural event coordinated by man-made malicious activities could take months to restore and cause significant national or international security concerns [145]. To perform advanced control functions, microgrids require advanced communication system devices that can be remotely monitored and controlled. However, advanced communication and control devices create vulnerabilities that could be exploited in malicious activities such as dial-up access to controllers, vendor supports, IT-controlled communication gear, access from corporate VPNs, database links, poorly configured firewalls, and peer utility links [146]. Therefore, more research in cyber-physical situational awareness and resilient control systems to guard against malicious activities for microgrids is needed.

**Standardization of resilience for microgrids:** Effective monitoring and control system designs that increase the resilience of microgrids against harmful events, require clarification and standardization of the definition of resilience, associated analytical methods, and metrics. In one definition, the term 'resilience' means the ability to prepare for and adapt to changing conditions and withstand and recover rapidly from disruptions. Resilience includes the ability to withstand and recover from deliberate attacks, accidents, or naturally occurring threats or incidents [147]. By that definition, together with additional terms defined in [148], a visualization of resilience measurements (resilience level) throughout all stages of power grids that are impacted by a cyber-physical event is illustrated in FIGURE 6. Although the concept of resilience for complex systems has been around for decades, in the power and energy community a universally accepted definition of power system resilience, metrics, and methodologies is not yet available [149]. And although there are multiple solutions that could be used to quantify system resilience, no single solution is applicable to all systems [150], [151]. Currently, the most common approach is to assess the overall social and economic impacts of a (potential) event on the system. Based on the monitoring of that system, control can be designed to minimize the identified impacts.

**Cyber-Social-Physical microgrids:** The future of system monitoring and control is moving toward the cyber-social-physical microgrid - resilient communities. An important trend in this field is the modeling and integrating the human behavior into the resilience control system, such as human behavior aware residential community EMS [152], customer-centric demand response, human-centric home EMS [153], human errors, resilient social network, etc. As human behavior and privacy have a large impact on the resilient objectives, these additional issues need to be considered in future monitoring and control research for microgrids.

### IV. CONCLUDING REMARKS

Microgrids are the most promising component of the power system capable of ensuring resilient energy services for critical infrastructure that is impacted by either natural disasters or man-made incidents. To guarantee real-time resilient operations of microgrids, monitoring and control functionalities are required. The trend in monitoring has recently shifted from normal situational awareness in forecasting, state estimation, and prediction to analysis of anomalies and detection of cyber-physical attacks. To respond to estimated or forecasted events, resilient control systems to ensure optimal power flow and guarantee system voltage and frequency stability were discussed. Although these monitoring and control systems have been investigated in power distribution and transmission systems, their applicability to and challenges for microgrids need to be addressed. Critically, there is not yet a widely accepted (consensus) definition, analytical methods, and associated metrics to consistently describe the resilience of power grids, especially for microgrids. This issue is tightly related to cyber-social-physical system design, monitoring, and control of microgrids. Therefore, future research should consider (1) cyber-physical situational awareness and resilient control systems to guard against malicious activities for microgrids, 2) the need for a standardized resilience framework for comprehensive and consistent resilience research in different monitoring and control layers of microgrids, (3) distinguishing and clarifying the definition of resilience for microgrids versus for distribution and transmission systems, and 4) human behavior in the system monitoring and control designs.

BIOGRAPHIES

**Tuyen V. Vu** is an Assistant Professor in the Department of Electrical and Computer Engineering at Clarkson University. He received his B.S. in electrical engineering from Hanoi University of Science Technology, Vietnam in 2012, and his Ph.D. in electrical engineering from Florida State University in 2016. From 2016 to 2018, he was a postdoctoral research associate and a research faculty at Florida State University - Center for Advanced Power Systems. Dr. Vu has worked on monitoring and control of cyber-physical microgrids for various applications since 2013.

**Bang L. H. Nguyen** received his B.S. and M.S. in electrical engineering from Ho Chi Minh University of Technology, Vietnam, in 2010 and 2013, respectively. In 2015, he was with Eastern International University, Binh Duong New City, Vietnam, as a lecturer. From 2016 to 2018, he was a research assistant in the power electronics and energy conversion lab (PEEC), Kyungpook National University, Korea. He is currently working towards his Ph.D. degree at Clarkson University. His areas of interest include the control and security of smart grids using deep learning.

**Zheyuan Cheng** received his B.S. degree in electrical engineering from Nanjing University of Aeronautics and Astronautics, China, in 2015. He is currently pursuing his Ph.D. degree in electrical and computer engineering at North Carolina State University, Raleigh. He joined the Advanced Diagnosis, Automation, and Control Laboratory at North Carolina State University in 2016. His research interests include distributed control systems, computational intelligence, and distributed optimization. He has been working on microgrid distributed resilient control related projects since 2016.

**Mo-Yuen Chow** is a Professor in the Department of Electrical and Computer Engineering at North Carolina State University. Dr. Chow's recent research focuses on distributed control and management on smart micro-grids, batteries, and mechatronics systems. Dr. Chow has established the Advanced Diagnosis, Automation, and Control Laboratory. He is an IEEE Fellow. He has received the IEEE Region-3 Joseph M. Biedenbach Outstanding Engineering Educator Award, the IEEE Industrial Electronics Society Anthony J Hornfeck Service Award. He is a Distinguished Lecturer of IEEE IES. Dr. Chow has been working on several projects related to the cyber-physical microgrids since 2008.

**Bin Zhang** is an Associate Professor in the Department of Electrical Engineering at the University of South Carolina, Columbia, SC. He received B.E. and M.E. degrees from Nanjing University of Science and Technology, Nanjing, China and Ph.D. degree from Nanyang Technological University, Singapore. Before he joined University of South Carolina, he was with General Motors, Detroit, MI, Impact Technologies, Rochester, NY, and Georgia Institute of Technology, Atlanta, GA. His research includes prognostics and health management, intelligent systems and controls. He has led and participated in some projects related to health management and control in recent years.